\documentclass[a4paper,reqno,12pt]{amsart}

\usepackage{bm}
\usepackage{amssymb}
\usepackage{tikz-cd}
\usepackage[utf8]{inputenc} 
\usepackage[T1]{fontenc}    
\usepackage{url}            
\usepackage{booktabs}       
\usepackage{amsfonts}       
\usepackage{nicefrac}       
\usepackage{microtype} 
\usepackage{enumerate}
\usepackage{lipsum}
\usepackage{calrsfs}
\usepackage{chngcntr}
\usepackage{pgfplots}
\usetikzlibrary{intersections}
\input xygraph
\input xy
\xyoption{all}
\RequirePackage{color}
\DeclareMathAlphabet{\pazocal}{OMS}{zplm}{m}{n}
\usepackage[all,color]{xy}
\usepackage{eqnarray,amsmath}

\def\o{\omega}
\def\a{\alpha}
\def\b{\beta}
\def\g{\gamma}
\def\d{\delta}
\def\l{\lambda}
\def\m{\mu}
\def\s{\sigma}
\def\E{\mathcal E}
\def\Z{\mathbb Z}
\def\N{\mathbb N}
\def\Q{\mathbb Q}
\def\R{\mathbb R}
\def\C{\mathbb C}
\def\K{\mathbb K}

\def\O{\mathbf{O}}
\newtheorem{lemma}{Lemma}
\newtheorem{corollary}{Corollary}
\newtheorem{definition}{Definition}
\newtheorem{theorem}{Theorem}
\newtheorem{proposition}{Proposition}

\title[2-dimensional perfect evolution algebras over domains]{Two dimensional perfect evolution algebras over domains}

\author[Y. Cabrera ]{Yolanda Cabrera Casado}
\address{Y. Cabrera Casado: Departamento de Matem\'atica Aplicada, E.T.S. Ingenier\'\i a Inform\'atica, Universidad de M\'alaga, Campus de Teatinos s/n. 29071 M\'alaga.   Spain. }
\email{yolandacc@uma.es}

\author[D. Mart\'{\i}n]{
 Dolores Mart\'{\i}n Barquero}
\address{D. Mart\'{\i}n Barquero: Departamento de Matem\'atica Aplicada. Escuela de Ingenier\'\i as Industriales. Universidad de M\'alaga, Campus de Teatinos. 29071 M\'alaga,   Spain.}
\email{dmartin@uma.es}
\author[C. Mart\'{\i}n]{Cándido Mart\'{\i}n Gonz\'alez*} 
\address{C. Mart\'{\i}n Gonz\'alez: Departamento de \'{A}lgebra, Geometr\'{i}a y Topolog\'{\i}a. Facultad de ciencias. Universidad de M\'{a}laga, Campus de Teatinos. 29010 Málaga, Spain}
\email{candido\_m@uma.es}

\makeatletter
\@namedef{subjclassname@2020}{%
  \textup{2020} Mathematics Subject Classification}
\makeatother

\thanks{ The  authors are supported by the Spanish Ministerio de Ciencia e Innovaci\'on through project  PID2019-104236GB-I00 and  by the Junta de Andaluc\'{\i}a  through projects  FQM-336 and UMA18-FEDERJA-119,  all of them with FEDER funds. 
}

\begin{document}
\maketitle

\begin{abstract}
We will study evolution algebras $A$ which are free modules of dimension $2$ over domains. Furthermore, we will assume that these algebras are perfect, that is $A^2=A$.
We start by making some general considerations about algebras over domains: they are sandwiched between a certain essential $D$-submodule and its scalar extension over the field of fractions of the domain. We introduce the notion of quasiperfect algebras and modify slightly the procedure to associate a graph  to an evolution algebra over a field given in \cite{ElduqueGraphs}.  Essentially, we introduce color in the connecting arrows, depending on a suitable criterion related to the squares of the natural  basis elements. Then we classify the algebras under scope parametrizing the isomorphic classes by convenient moduli. 
\end{abstract}

\section{Introduction}
There is a large number of publications  studying 2-dimensional evolution algebras (\cite{CLOR}, \cite{RozikovMurodov}, \cite{Murodov}, \cite{CelorrioVelasco}, \cite{Ahmed}, \cite{Ceballos}, \cite{CGMMS}, etc.), among them we could highlight those works that deal with their classification. So,  in the works of \cite{Murodov} and \cite{CLOR} the   evolution algebras of dimension two over the reals and complex respectively are classified,  the paper \cite{CelorrioVelasco} contains the  classification of the evolution algebras of dimension two and three over the field $\K$ , where $\K=\R$ or $\C$. The work \cite{Ahmed} is addressed to the case in which the ground field $\K$ is algebraically closed. In \cite{CGMMS} the classification of the evolution algebras of dimension  two over arbitrary fields is provided. We can find more information about the evolution of the research in the field of Tiam's evolution algebras in \cite{ceballo}.  We are contributing our grain of sand by studying the perfect two-dimensional  evolution algebras over domains. So, in this paper, the word domain will stand for a commutative ring such that $xy=0$ implies $x=0$ or $y=0$. If $D$ is a domain, an \emph{evolution algebra} over $D$ is a $D$-algebra which is free as $D$-module and has a basis such that the product of any two different elements is zero. Such basis is called a \emph{natural basis}. Of course, any one-dimensional $D$-algebra is an evolution algebra.  Evolution algebras over domains are much more involved than evolution algebras over fields.
For instance, there are two one-dimensional evolution algebras over fields up to isomorphism: the ground field with zero multiplication and the ground field with its usual multiplication. But if you consider a domain $D$ and
a one-dimensional evolution $D$-algebra there are more isomorphic classes of evolution algebras.
Ruling out the trivial one, we have define in $D$ a multiplication with $1^2=d$ where $d\in D^*:=D\setminus\{0\}$. Denote this algebra by $D_d$. The product in $D_d$ is $x\cdot y=xyd$ for any $x,y\in D$.
If $f\colon D_d\to D_e$ is an isomorphism, it is a $D$-module isomorphism hence $f(a)=af(1)$ for any $a$.
So $f$ is the multiplication times an element $x:=f(1)$. Consequently $x\in D^\times$ (the group of invertibles of $D$). 
But also $f(uvd)=f(u)f(v)e$ whence $f(d)=f(1)^2e$ or $d=xe$.  
Thus the isomorphic condition is $D_d\cong D_e$ if and only if there exists $x\in D^\times$ such that $d=xe$.
So the isomorphic classes of nontrivial one dimensional evolution algebras over domains with nonzero product
is in one to one correspondence with the set $D^*/D^\times$, that is, the
set of equivalence classes of $D^*$ modulo the action $D^\times\times D^*\to D^*$ such that $x\cdot d=x d$ for any $x\in D^\times$, $d\in D^*$.
For instance, if $D=\Z$ we have $\Z^*/\{\pm 1\}\cong\N^*=\{1,2,\ldots\}$.
Thus there are countably many isomorphic classes of one-dimensional evolution algebras over certain domains. If we consider the domain $D:=\K[x]$ of polynomials in one
indeterminate over the field $\K$ then the isomorphic classes of one-dimensional
evolution $D$-algebras are in one to one correspondence with the set of monic
polynomials of $\K[x]$.

One of the elements of this study is the use of moduli sets from certain classes of algebras: the idea is to parametrize the algebras of a class by tuples of parameters ranging in a given space. It turns out that in some cases the tuples range in curves or surfaces or other varieties. The different algebras in the same isomorphic class may happen to be in a curve of an affine plane and the different curves fill the space modulo the restrictions on the parameters imposed by the class of algebras. This may be seen as a bundle in the category of sets.

This paper is organized as follows. In Section \ref{sec:headings} we prove some results on algebras over domains in Proposition \ref{leurc} and  Corollary \ref{reed}. We introduce the class of quasiperfect algebras and give necessary and sufficient conditions for an evolution algebra to be perfect or quasiperfect in terms of natural bases (in Proposition \ref{basischange}). Next, we define  the required terminology to be able to classify our $D$-algebras using moduli sets.   Finally, in Section \ref{perfectcase} we associate a colored directed graph to a quasiperfect evolution $D$-algebra and give the classification theorem  of perfect two-dimensional evolution algebras over a domain (Theorem \ref{teoremon}).

\smallskip

\section{Preliminaries and previous results}
\label{sec:headings}

If $D$ is a (commutative) domain we will denote by $Q:=Q(D)$ the field of fractions of $D$. For any $D$-module $M$ we will construct $M_Q$ the $Q$-module of fractions $M_Q:=\{\frac{m}{d}\colon m\in M, d\in D\setminus\{0\}\}$ with the usual operations of sum, product and product by elements of $Q$. To be more specific we consider the set of couples $(m,d)$ with $m\in M$, $0\ne d\in D$ modulo the equivalence relation $$(m,d)\equiv(m',d') \text{ iff } \exists t\in D\setminus\{0\}, t(d'm-dm')=0.$$
Then we denote the equivalence class of $(m,d)$ in the usual way: $\frac{m}{d}$ and $M_Q$ is a $Q$-vector space relative to the usual sum of fractions and $\frac{d_1}{d_2}\frac{m}{d_3}:=\frac{d_1m}{d_2d_3}$. If $M$ is a torsion-free $D$-module, there is a canonical monomorphism of $D$-modules
$M\to M_Q$ such that $m\mapsto\frac{m}{1}$. Usually we will denote $\frac{m}{1}=m$ so that the elements of $M$ and their images in $M_Q$ will be identified. Unless otherwise stated, we will work allthrough this paper with torsion-free $D$-modules.

A well known property of $M$ is that a set of vectors $\{e_i\}\subset M$ is linearly independent if and only if its image in $M_Q$ is $D$-linearly independent. As well if $\{u_i\}_{i\in I}$ is linearly independent in the $Q$-vector space $M_Q$ and $u_i=\frac{m_i}{d_i}$ for any $i\in I$, then the set 
of numerators $\{m_i\}_{i\in I}$ is linearly independent in $M$. We also have:

\begin{proposition} \label{leurc}
Let $D$ be a domain, $Q$ its field of fractions, $M$ a torsion-free $D$-module and $M_Q$ the $Q$-module of fractions of $M$.
\begin{enumerate}
\item[(a)] If $\{u_i\}_{i\in I}$ is a basis of $M_Q$ (as a $Q$-vector space)  and $u_i=\frac{m_i}{d_i}$ for $i\in I$, then 
$\{m_i\}$ is also a basis of $M_Q$ and 
a maximal linearly independent subset of $M$.
\item[(b)] A set $\{m_i\}\subset M$ is a maximal linearly independent subset of $M$ if and only if  $\{m_i\}$ is a basis of the $Q$-vector space $M_Q$. 
\item[(c)] If $M$ is a $D$-module and
$\{m_i\}$ a maximal linearly independent subset of $M$, then $\oplus_i Dm_i$ is an essential submodule of $M$.
\item[(d)] Assume that $M$ is a free $D$-module with a finite basis $\{e_i\}_{i\in I}$. A maximal linearly independent subset $\{u_i\}_{i\in I}$ of $M$ is a basis of $M$ if and only if the determinant of the change of basis matrix is in $D^\times$.
\end{enumerate}
\end{proposition}
Proof. For the first assertion, take into account that the set of the $m_i$'s is linearly independent also in $M_Q$ because if we have 
$\sum \l_i m_i=0$ then we can write $\sum\l_i d_i \frac{m_i}{d_i}=\sum\l_i d_i u_i=0$ hence for any $i$ we have $\l_id_i=0$ which implies $\l_i=0$. The set $\{m_i\}$ is also a system of generators of the $Q$-vector space $M_Q$ because any $x\in M_Q$
can be written as 
$$x=\sum\l_iu_i=\sum\l_i\frac{m_i}{d_i}=
\sum\l_i\frac{1}{d_i} m_i.$$
Let us prove now that $\{m_i\}$ is a maximal among the linearly independent subsets of $M$: if $\{m_i\}\subsetneq T$ for a linearly independent subset $T\subset M$, then $T$ contains properly a basis of $M_Q$ and is linearly independent, a contradiction.

For the second assertion, take  a maximal linearly independent set $\{m_i\}\subset M$.
We know that $\{m_i\}$ is linearly independent also in $M_Q$. If this set is not a basis, there is some $x\in M_Q$ such that 
$\{m_i\}\cup\{x\}$ is again linearly independent. If $x=\frac{z}{d}$ then 
$\{m_i\}\cup\{z\}$ is a linearly independent subset of $M$ contradicting the maximality of $\{m_i\}$. Reciprocally, if $\{m_i\}\subset M$ is a basis of $M_Q$ we know that $\{m_i\}$ is linearly independent in $M$. To prove the third assertion take a nonzero submodule $N$ of $M$. We must prove that $N\cap (\oplus_i Dm_i)\ne 0$. Take $0\ne n\in N$. Since $\{m_i\}$ is a basis of $M_Q$ we have 
$dn=\sum_i d_i m_i$ for some $d,d_i\in D$ (and $d\ne 0$). Thus $0\ne dn\in N\cap(\oplus_i Dm_i)$.
Finally we prove the fourth assertion. We have
$u_i=a_i^je_j$ (using Einstein summation  convention) for any $i\in I$. If $\det[(a_i^j)]\in D^\times$, then there are scalars $b_i^j\in D$ such that
$e_i=b_i^je_j$ for any $i$. Whence $\{u_i\}_{i\in I}$ is a basis of $M$. Reciprocally, if $\{u_i\}_{i\in I}$ turns out to be a basis we may write $e_i=b_i^je_j$
for suitable scalars $b_i^j\in D$. But then the matrices $(a_i^j)$ and $(b_i^j)$ have product $1$, that is, $(a_i^j)(b_i^j)=1$. This implies that the determinant of each such matrix is an invertible element of $D$. $\square$
\medskip 

\begin{corollary}\label{reed}
Let $A$ be a $D$-algebra, then there is a maximal linearly independent subset $\{a_i\}$ of $A$ (in fact a basis of the $Q$-vector space $A_Q$) such that $A$ is contained as $D$-module in a sandwich $$\oplus D a_i\subset A\subset A_Q,$$
and $\oplus Da_i$ is an essential $D$-module of $A$.
\end{corollary}
Proof.
Take a basis $\{u_i\}$ of $A_Q$ as a $Q$-vector space. If $u_i=a_i/d_i$ then $\{a_i\}$ is a maximal linearly independent subset of $A$ by Proposition \ref{leurc}. Also $\oplus Da_i$ is essential as a $D$-submodule of $A$ by Lemma \ref{leurc}. Note that $A$ is a torsion-free $D$-module: if $da=0$ for some nonzero $d\in D$, write $a=\sum_i q_i u_i$ as a linear combination of the $u_i$'s. Then we have $0=d\sum_i q_i u_i$ whence $dq_i=0$ for any $i$. Since $d\ne 0$ we have $q_i=0.$
$\square$
\medskip

In the situation above $A\subset A_Q$, if $A$ is perfect ($A^2=A$), then $A_Q$ is also perfect. However, we may have $A_Q$ perfect and $A$ not.
For instance, consider $A=2\Z\times\Z$ with componentwise multiplication. Then $A_\Q\cong\Q\times\Q$ which is perfect but $A$ is not. 

An example of the situation described in Corollary \ref{reed} is given by taking $D=\Z$ and $A=\{(\frac{x}{2},y)\colon x,y\in\Z\}=\Z(\frac{1}{2},0)\oplus\Z(0,1)$. Then $A$ is a two-dimensional free $\Z$-module and $A_\Q=\Q(1,0)\oplus\Q(0,1)$ being 
$$\Z(1,0)\oplus\Z(0,1)\subsetneq A\subsetneq\Q(1,0)\oplus\Q(0,1).$$

\begin{definition}\rm 
An algebra $A$ over a domain $D$ will be termed {\em quasiperfect} if $A_Q^2=A_Q$.
\end{definition}

Let $D$ be a domain and $\E$ a free $D$-module  such that $\E$ is an evolution algebra with natural basis $\{e_1,e_2\}$.
We will need the following proposition whose first item is exactly the same as in the case of evolution algebras over fields:
\begin{proposition}\label{basischange}
Assume that $\E$ is an evolution algebra over a domain $D$ with a finite natural basis $\{e_i\}_{i\in I}$. Let $\o_i^j\in D$ be the structure constants, that is, $e_i^2=\o_i^j e_j$ (using Einstein summation  convention). Then we have: 
\begin{enumerate}
\item $\E^2=\E$ if and only if the matrix $(\o_i^j)$ is invertible.
Moreover, for any other natural basis $\{f_i\}_{i\in I}$ there is a permutation $\sigma$ of $I$ such that $f_i=k_ie_{\sigma(i)}$ and each $k_i\in D^\times$.
\item $\E$ is quasiperfect if and only if the determinant of $(\o_i^j)$ is nonzero. As in the previous case, for any other natural basis $\{f_i\}_{i\in I}$ there is a permutation $\sigma$ of $I$ such that $f_i=k_ie_{\sigma(i)}$ and each $k_i\in D^\times$.
\end{enumerate}
\end{proposition}
Proof. From $\E=\E^2$ we deduce that $e_i=x_i^je_j^2$ for any $i\in I$. Then $e_i=x_i^j\o_j^k e_k$ so that $x_i^j\o_j^k=\delta_i^k$ (Kronecker Delta). Thus the matrix $(\o_i^j)$ is invertible. 
Reciprocally if $(\o_i^j)$ is invertible
the linear map such that $e_i\mapsto e_i^2$ is
an isomorphism whence $\E^2=\E$. Assuming the perfection of $\E$, if $\{f_i\}$ is another natural basis and we write $f_i=a_i^j e_j$, then
for $i\ne j$ we have $0=f_if_j=a_i^ke_k a_j^q e_q=a_i^ka_j^q e_ke_q=a_i^ka_j^q \delta_{kq}\o_k^s e_s=a_i^ka_j^k \o_k^s e_s$ whence $a_i^ka_j^k\o_k^s=0$ for any $s$ and any couple $(i,j)$ with $i\ne j$.
Since the matrix $(\o_i^j)$ is invertible, we consider its inverse matrix $(\tilde \o_i^j)$, so we have 
$\o_i^j\tilde \o_j^k=\delta_i^k$. Then, from 
$a_i^ka_j^k\o_k^s=0$ we get 
$a_i^ka_j^k\o_k^s\tilde\o_s^q=0$ or 
$a_i^ka_j^k\delta_k^q=0$. Thus $a_i^qa_j^q=0$ for any $q$ provided $i\ne j$. So in each column and each row of the matrix $(a_i^j)$ there is a unique nonzero element.
Consequently, $f_i=k_i e_{\sigma(i)}$ for a certain permutation $\sigma$ of $I$. Now, the coefficients $k_i$ are invertible in $D$ since the determinant of the matrix
of basis change is invertible. 
Let us prove now the second assertion. If we have $\E_Q^2=\E_Q$ hence the matrix $(\o_i^j)$ is invertible in $Q$. Its determinant is a nonzero element of $D$. Reciprocally if $\det(\o_i^j)\ne 0$ then it is invertible in $Q$ so that $\E$ is quasiperfect. In this case,
if $\{e_i\}$ and $\{f_i\}$ are natural bases of $\E$, then there is a permutation $\sigma$ and nonzero elements $k_i\in Q$ such that $f_i=k_ie_{\s(i)}$ for any $i$. Also $f_i=a_i^je_j$ for certain $a_i^j\in D$ hence $a_i^{\s(i)}=k_i$ and $a_i^j=0$ if $j\ne\s(i)$.
In any case $k_i\in D^\times$.
$\square$
\begin{corollary}
Any perfect evolution algebra over a domain is quasiperfect.
\end{corollary}
 Next, we prove a result which generalizes the situation of perfect evolution algebras over fields.
\begin{lemma}\label{veraiv}
Assume that $\E$ is a quasiperfect evolution algebra over a domain $D$ with a finite natural basis $\{e_i\}_{i\in I}$.
Then the following numbers do not depend on the natural basis chosen:
\begin{enumerate}
    \item The number of nonzero entries in the structure matrix
    $(\o_i^j)$.
    \item The number of nonzero entries in the diagonal of $(\o_i^j)$.
    \item The number of invertible elements in $(\o_i^j)$.
    \item The number of invertible elements in the diagonal of $(\o_i^j)$.
\end{enumerate}
\end{lemma}
Proof. For our original natural basis  we have $e_i^2=\o_i^je_j$.
Take now any other natural basis $\{f_i\}$. There is a permutation $\sigma$ of $I$ such that 
$f_i=k_i e_{\sigma(i)}$ for some
invertible elements $k_i\in D^\times$. It follows that 
$f_i^2=\tau_i^q f_q$ where 
$$\tau_i^q=\frac{k_i^2}{k_q}\o_{\sigma(i)}^{\sigma(q)}.$$
So, the number of nonzero (respectively invertible) elements in the matrix $(\o_i^j)$ coincides with the  number of nonzero (respectively invertible) elements in the matrix $(\tau_i^j)$. Similarly for the diagonal elements. $\square$
\medskip

We can associate a colored graph to any quasiperfect evolution algebra over a domain. The key idea is of course that of \cite{ElduqueGraphs}  but slightly modified. We fix a natural basis $\{e_i\}$ of an evolution $D$-algebra $\E$. Let $\omega_i^j$ be the structure constants. Then, we consider the graph whose vertices are in bijection with $\{e_i\}$ and we draw a black edge from vertex
$i$ to vertex $j$ if $\omega_i^ j\in D^\times$. We draw a blue edge from vertex $i$ to vertex $j$ if $\omega_i^j\in D^*\setminus D^\times$. The (isomorphic class) of the associated graph does not depend on the chosen natural basis because of Lemma \ref{veraiv}. So, the graph associated to an evolution algebra over a domain is a colored directed graph.

\begin{proposition}\label{shell} Let $D$ be a domain and $Q$ its field of fractions.
If $\E$ is a quasiperfect evolution algebra over $D$ of finite dimension and $\{m_i\}$ is maximal linearly independent subset of $\E$ such that $m_im_j=0$ for $i\ne j$, then $\{m_i\}$ is a natural basis of $\E$.
\end{proposition}
Proof.
Fix a natural basis $\{e_i\}$ of $\E$.
We know that $\{m_i\}$ is a basis of $\E_Q$ by Proposition \ref{leurc} and it is a natural basis. Then $m_i=k_ie_{\s(i)}$ for some $k_i\in Q^\times$ and a permutation $\s$. On the other hand, there are expressions 
$m_i=a_i^j e_j$ where $a_i^j\in D$. Consequently $a_i^{\s(i)}=k_i$ so that each $k_i\in D^\times$. Furthermore, $a_i^j=0$ if $j\ne\s(i)$. Thus, $\{m_i\}$ is a natural basis of $\E$. $\square$\medskip

\section{The moduli set.}
We have to introduce some terminology in order to present what we will call the moduli set for the different classes of algebras. 
Given a class $C$ of algebras  we will say that a set $S$ is a \emph{moduli} for $C$ (or a \emph{moduli set}) if there is a one-to-one correspondence between the isomorphic classes of algebras of $C$ and the elements of $S$. In some occasions the moduli of a class will have an additional algebraic or geometric structure. In the category of sets, we recall that a \emph{bundle} is an epimorphism $\pi\colon E\to B$ where 
$E$ is called \emph{total set}, $B$ is the \emph{base set}  and for every $b\in B$, $p^{-1}(b)$ is the \emph{fiber} over $b$. 
By a \emph{cross section} of $\pi$ we understand a right inverse $s\colon B\to E$. 
We have an example of bundle if we consider as the total set a class of algebras $C$ and as the base set the set $C/\cong$ of isomorphic classes of the algebras in $C$. Then, $\pi$ maps any algebra to its isomorphic class. The fiber of an element represents an isomorphic class of algebras. If we specify a cross section of the bundle, then this is equivalent to give a representative of the isomorphic class of any element of $C$. Thus, the classification problem of $C$ under isomorphism consists just in giving a cross section of the corresponding bundle. The moduli set of the classification is the base set of the corresponding bundle.  

\subsection{Direct limits} Some of the cases of our classification are based on direct limits. If $M$ is an abelian monoid and $n\in\Z$ we consider the direct system 
$$M\buildrel{(\cdot)^n}\over\to M\buildrel{(\cdot)^n}\over\to\cdots\to M\buildrel{(\cdot)^n}\over\to\cdots$$ where $(\cdot)^n$ is the homomorphism such that $g\mapsto g^n$,  then we will denote
the direct limit of such system by $\displaystyle\lim_{\to n}M$. Recall that this monoid can be described as follows: 
consider the sequence of monoids $\{M_i\}_{i\in\N}$ such that $M_i:=M$ for any $i$. Then, in the disjoint union 
$\sqcup_i M_i$ we define an equivalence relation: if $x\in M_i$ and $y\in M_j$ we say that $x\sim y$ if and only if 
$y^{n^k}=x^{n^h}$ for some naturals $k,h$. We denote the equivalence class of $g\in M$ by $[g]$.
So $\lim_{\to n} M$ is the quotient of $\sqcup_i M_i$ modulo $\sim$.
A particular case of this arises if we take $M:=D^\times/(D^\times)^{[q]}$ where 
$(D^\times)^{[q]}:=\{x^q\colon x\in D^\times\}$ with $q \in \N^*$.
This specific $M$ is a group and its elements
are  equivalence classes $\bar\l$ with $\l\in D^\times$. We have $\bar\l= \bar\m$ if and only if
$\l=\m\ r^q$ for some $r\in D^\times$.
 We can consider 
$\lim_{\to n}M$ whose elements are the equivalence classes $[\bar\l]$ with $\bar\l\in M$. 

We have a canonical group homomorphism $D^\times\to \lim_{\to n}M$ such that $\l\to[\bar\l]$. Two elements $\l,\m\in D^\times$ are said to have the same image in $\lim_{\to n}M$ if $[\bar\l]=[\bar\m]$. For instance, consider the group $\lim_{\to 2}D^\times/(D^\times)^{[3]}$.
Since $2$ and $3$ are coprime the transition homomorphisms $(\cdot)^2$ are isomorphisms. Then $\l$ and $\m$ have the same image in $\lim_{\to 2}D^\times/(D^\times)^{[3]}$ if $\bar\l^{2^\a}=\bar\m^{2^\b}$ for some naturals $\a,\b\ge 0$. If $\a=\b$, then $\bar\l=\bar\mu$. If $\a>\b$, then $\bar\m=(\bar\l)^{2^{\a-\b}}$ and if $\a<\b$
we have $\bar\l=(\bar\m)^{2^{\b-\a}}$. In each case, replacing $2^{\a-\b}$ or $2^{\b-\a}$ with its remainder $\xi$ modulo $3$ we have 
$\bar\m=\bar\l^\xi$ or 
$\bar\l=\bar\m^\xi$ where $\xi=1,2$. So, 
there is some $r\in D^\times$ such that $\m=\l r^3$ or $\m=\l^2 r^3$ (observe that the other possibility $\l=\m^2 s^3$ is a consequence of $\m=\l^2 r^3$). Summarizing we have the following lemma.
\begin{lemma}
The elements $\l,\m\in D^\times$ have the same image in the group $\lim_{\to 2}D^\times/(D^\times)^{[3]}$ if and only if $\m=\l r^3$ or $\m=\l^2 r^3$ for some
$r\in D^\times$.
\end{lemma}

In this case, we have a bundle $\pi\colon D^\times\to \lim_{\to 2}D^\times/(D^\times)^{[3]}$ such that 
$\l\mapsto [\bar\l]$. The fiber of some $[\bar\l]$ is the set of all $\m\in D^\times$ such that $\m=\l r^3$ or $\m=\l^2 r^3$ for some $r\in D^\times$.
\subsection{Algebraic sets}
Assume that we have an action $D^\times\times X\to X$ where
$X$ is some subset of $D^n$. So, for $t\in D^\times$
and $a=(a_1,\ldots,a_n)\in X$ we might have 
$t\cdot(a_1,\ldots,a_n)=(x_1,\ldots,x_n)$ where each $x_i$ is a polynomial in $t$ with coefficients
in $D$. Consequently, we may write $x_i=p_i(t)$ where $p_i$ is the mentioned polynomial. The orbit of $a$ under the above action
is in a curve $x_i=p_i(t)$. More precisely, 
the orbit is contained in the image of  the map
$c\colon Q\to Q^n$ such that $t\mapsto (p_1(t),\ldots,p_n(t))$ where $Q$ is the field of fractions of $D$. Then the Zariski closure of the image of  $c$ is an algebraic set $V\subset Q^n$ and the orbit of $a\in X$ is just $V\cap X$. This setting will appear in our classification of evolution algebras.\medskip
\begin{enumerate}[(I)]\label{tnerap}
\item There is a family $C$ of two-dimensional evolution algebras depending on two parameters, namely those algebras $A(1,\l,\m)=D\times D$ where $D$ is a fixed domain and the product in the algebra being $e_1^2=e_2$, $e_2^2=\l e_1+\m e_2$ with $\l\in D^\times$ and $\m\in D^*\setminus D^\times$. As we will see in the next section the isomorphism conditions for algebras of $C$ is $A(1,\l,\m)\cong A(1,\l',\m')$ if and only if there is some $k\in D^\times$ such that 
$\l'=k^3\l$ and $\m'=k^2\m$. This induces an action $D^\times\times X_0\to X_0$ where $X_0=D^\times\times(D^*\setminus D^\times)$ given by $k\cdot(\l,\m)=(k^3\l,k^2\m)$. 
Note that the isotropy subgroup of any $(\l,\m)\in X_0$ is trivial. This implies that the cardinal of each orbit agrees with that of $D^\times$.

For a fixed $(\l,\m)\in X_0$ we can consider $c\colon Q\to Q^2$ such that
$c(k)=(k^3\l,k^2\m)$. The Zariski closure of the image of $c$ is $V(I)$, the algebraic set of zeros of the ideal $I\triangleleft Q[x,y]$ generated by the polynomial $\m^3x^2-\l^2y^3$. So, it is a curve $c_{\l,\m}$ of $Q^2$. Thus the cardinal of the orbit of $(\l,\m)$
agrees with that of the set of points of $c_{\l,\m}$
lying on $X_0$. Any point of $X_0$ is in some curve $c_{\l,\m}$, in fact, $(\l,\m)\in c_{\l,\m}$.
Denote by $c_{\l,\m}^*$ the section $c_{\l,\m}^*:=c_{\l,\m}\cap X_0$. Then
 if $(\l',\m')\notin c_{\l,\m}$ we have 
$c_{\l,\m}^*\cap c_{\l',\m'}^*=\emptyset$. So, $X_0$ is the disjoint union of all sections $c_{\l,\m}^*$ and we have a bundle (in the category of sets) $p\colon X_0\to X_0/D^\times$ in which $p(\l,\m)=\hbox{\rm orb}(\l,\m)$ can be identified with  $c_{\l,\m}^*$. The fibers of this bundle are the points in one specific curve so the fibers represent classes of isomorphic algebras. 

\begin{lemma}\label{upac}
The cardinal of each orbit of $X_0/D^\times$ is $\vert D^\times\vert$ and agrees with that of the set of points of
the curve $c_{\l,\m}\equiv \m^3x^2-\l^2y^3=0$ in $D^\times\times(D^*\setminus D^\times)$.
The orbit set $X_0/D^\times$ is the base space of a bundle $p\colon X_0\to X_0/D^\times$ where the fibers represent classes of isomorphic algebras. So $X$ is a disjoint union of sections $c_{\l,\m}\cap X_0$.
\end{lemma}

The previous bundle can be \lq\lq lifted\rq\rq\ to specific fields, for instance, in the real case 
we may consider the plane with the axis removed: $\Gamma:=\R^\times\times\R^\times$. Denote by $\O$ the origin $\O=(0,0)$. Consider also the curves 
$c_{\l,\m}$ each one of which, is the zero set of $\m^3x^2-\l^2 y^3$.
Then $\Gamma=\sqcup_{\l,\m} c^*_{\l,\m}$ is the disjoint union of the perforated curves $c_{\l,\m}^*:=c_{\l,\m}\tiny\setminus\{\O\}$ with $\l,\m\ne 0$. Each such curve $c_{\l,\m}$ cuts the line $x=1$ in an unique point: $(1, \mu/\root 3\of {\l^2})$. Then 
we consider $\pi\colon\R^\times\times\R^\times \to \R^\times$
where $\pi(\l,\m)$ can be defined as the intersection of $c_{\l,\m}$ with the vertical line $x=1$. In other words $\pi(\l,\m)=\frac{\m}{\root 3\of {\l^2}}$. Since for any nonzero $t$ we have $\pi(1,t)=t$ the map $\pi$ is surjective and defining the curve $\tilde c_t:=c_{\l,\m}$ for any $(\l,\m)\in\pi^{-1}(t)$ we have $$\R^\times\times\R^\times=\displaystyle\bigsqcup_{t\in \R^\times}\tilde c_t^*$$ where $\tilde c_t^*:=\tilde c_t\setminus\{0\}$. We could paraphrase this by saying that the axis-less plane is a disjoint union (indexed in $\R^\times$) of perforated curves. 

\[
\begin{matrix}
\begin{tikzpicture}[scale=0.4]
\newcommand\mier{2}
\begin{axis}[
   		xmin=-1, xmax=1,
   		ymin=-1, ymax=1,
   		xtick distance=1, ytick distance=4 ]
\addplot [domain=-1:1, samples=100, name path=f, thick, color=black!50]{\mier*0.1*exp(0.6*ln(abs(x)))};

\end{axis}
\end{tikzpicture} & 
\begin{tikzpicture}[scale=0.4]
\newcommand\mier{4}
\begin{axis}[
   		xmin=-1, xmax=1,
   		ymin=-1, ymax=1,
   		xtick distance=1, ytick distance=4 ]
\addplot [domain=-1:1, samples=100, name path=f, thick, color=black!50]{\mier*0.1*exp(0.6*ln(abs(x)))};

\end{axis}
\end{tikzpicture} &
\begin{tikzpicture}[scale=0.4]
\newcommand\mier{4}
\begin{axis}[
   		xmin=-1, xmax=1,
   		ymin=-1, ymax=1,
   		xtick distance=1, ytick distance=4 ]
\addplot [domain=-1:1, samples=100, name path=f, thick, color=black!50]{\mier*0.1*exp(0.6*ln(abs(x)))};

\end{axis}
\end{tikzpicture}\\
\begin{tikzpicture}[scale=0.4]
\newcommand\mier{6}
\begin{axis}[
   		xmin=-1, xmax=1,
   		ymin=-1, ymax=1,
   		xtick distance=1, ytick distance=4 ]
\addplot [domain=-1:1, samples=100, name path=f, thick, color=black!50]{\mier*0.1*exp(0.6*ln(abs(x)))};
\end{axis}
\end{tikzpicture} & 
\begin{tikzpicture}[scale=0.4]
\newcommand\mier{8}
\begin{axis}[
   		xmin=-1, xmax=1,
  		ymin=-1, ymax=1,
   		xtick distance=1, ytick distance=4 ]
\addplot [domain=-1:1, samples=100, name path=f, thick, color=black!50]{{\mier*0.1*exp(0.6*ln(abs(x)))}};
\end{axis}
\end{tikzpicture} &
\begin{tikzpicture}[scale=0.4]
\newcommand\mier{10}
\begin{axis}[
   		xmin=-1, xmax=1,
   		ymin=-1, ymax=1,
  		xtick distance=1, ytick distance=4 ]
\addplot [domain=-1:1, samples=100, name path=f, thick, color=black!50]{\mier*0.1*exp(0.6*ln(abs(x)))};
\end{axis}
\end{tikzpicture}\\
\begin{tikzpicture}[scale=0.4]
\newcommand\mier{12}
\begin{axis}[
   		xmin=-1, xmax=1,
   		ymin=-1, ymax=1,
   		xtick distance=1, ytick distance=4 ]
\addplot [domain=-1:1, samples=100, name path=f, thick, color=black!50]{\mier*0.1*exp(0.6*ln(abs(x)))};
\end{axis}
\end{tikzpicture} & 
\begin{tikzpicture}[scale=0.4]
\newcommand\mier{14}
\begin{axis}[
   		xmin=-1, xmax=1,
   		ymin=-1, ymax=1,
   		xtick distance=1, ytick distance=4 ]
\addplot [domain=-1:1, samples=100, name path=f, thick, color=black!50]{\mier*0.1*exp(0.6*ln(abs(x)))};
\end{axis}
\end{tikzpicture} &
\begin{tikzpicture}[scale=0.4]
\newcommand\mier{16}
\begin{axis}[
   		xmin=-1, xmax=1,
   		ymin=-1, ymax=1,
   		xtick distance=1, ytick distance=4 ]
\addplot [domain=-1:1, samples=100, name path=f, thick, color=black!50]{\mier*0.1*exp(0.6*ln(abs(x)))};
\end{axis}
\end{tikzpicture}
\end{matrix}
\]
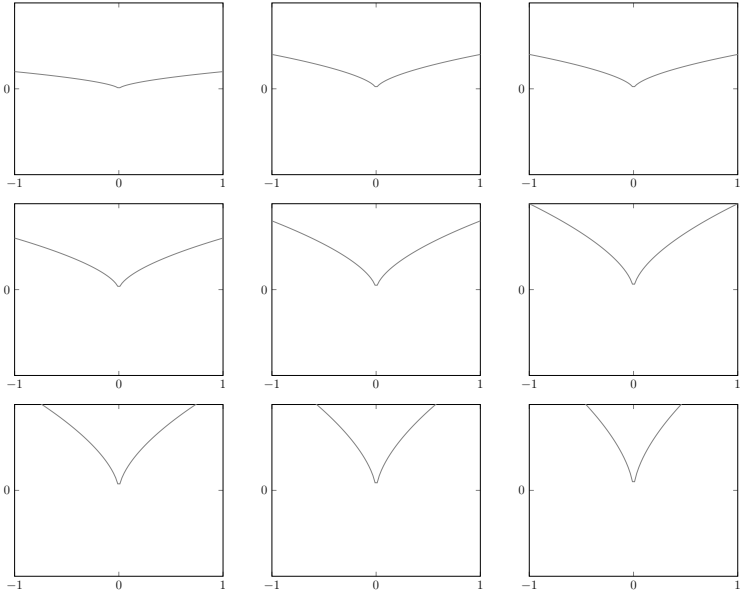
\centerline{Some of the curves $c_{\l,\m}$.}
\newline 

\begin{figure}
\[
 \begin{tikzpicture}

  \begin{axis}[
  		xmin=-3, xmax=3,
  		ymin=-10, ymax=10,
  		xtick distance=1, ytick distance=4 ]

\foreach \i in {1,...,25}
{
    \addplot [domain=-2.5:2.5, samples=100, name path=f, thick, color=black!50]
        {\i*exp(0.6*ln(abs(x)))};
        \addplot [domain=-2.5:2.5, samples=100, name path=f, thick, color=black!50]
       {-\i*exp(0.6*ln(abs(x)))};
}
 
   \draw [opacity=0.4] (axis cs:{1,-10}) -- (axis cs:{1,10});
    \draw [dashed, opacity=0.4] (axis cs:{-2.5,0}) -- (axis cs:{2.5,0});
     \draw [dashed, opacity=0.4] (axis cs:{0,-8}) -- (axis cs:{0,8});

   \end{axis}

  \end{tikzpicture}
\]
\caption{The family $c_{\l,\m}$.}
\end{figure}

If we consider for instance the domain $D=\Z[\sqrt{3}]\subset\R$, we know that $D^\times$ consists of all $a+b\sqrt{3}\in D$ such that $a^2-3b^2=1$ which has infinite cardinal. Thus, in each orbit of
$X_0=\Z[\sqrt{3}]^\times\times(\Z[\sqrt{3}]^*\setminus \Z[\sqrt{3}]^\times)$ under the action of $\Z[\sqrt{3}]^\times$ there are infinitely many elements. The orbits of this actions are in one-to-one correspondence with the points of the real line of the form $\frac{\m}{\root 3\of{\l^2}}$ where
$\m\in \Z[\sqrt{3}]^*$ and $\l\in\Z[\sqrt{3}]^\times$.
So, we can consider the monoid $\mathcal{M}:=\{\frac{\m}{\root 3\of {\l^2}}\colon \m\in\Z[\sqrt{3}]^*, \l\in \Z[\sqrt{3}]^\times\}\subset\bar\Q$ and the orbits of the action $\Z[\sqrt{3}]^\times\times X_0\to X_0$ are in one-to-one correspondence with the monoid $\mathcal{M}$. Consequently, the isomorphic classes of algebras of type $A(1,\l,\m)$ are in one-to-one correspondence with $\mathcal{M}$.

To manage the general case of Lemma \ref{upac} we have to make the following considerations. 
Pick a domain $D$, let $Q$ be its field of fractions and $\bar Q$ the algebraic closure of $Q$. Take the multiplicative monoid 
$\mathcal{M}:=\{
\frac{\m}{k}\colon \m\in D^*, 
\exists\l\in D^\times,\ k^3=\l^2
\}\subset\bar Q$. 
Then, we consider $\m_3(\mathcal{M}):=
\{x\in\mathcal{M}\colon x^3=1\}$ 
and define the quotient monoid
\begin{equation}\label{sapr1}
\bar{\mathcal{M}}:=\mathcal{M}/\m_3(\mathcal{M})
\end{equation}
whose elements are equivalence classes $[\l]$. 
It is interesting to note that $\m_3(\mathcal{M})=\m_3(\bar Q)$: indeed, if
$x\in\m_3(\bar Q)$ we have $x=\frac{1}{x^2}\in\mathcal{M}$ since $(x^2)^3=1^2$.
Furthermore,  there is a bijection $X_0/D^\times\cong\bar{\mathcal{M}}$ such that the equivalence class of $(\l,\m)$ in $X_0$ maps to $[\frac{\m}{\root 3\of{\l^2}}]$. Definitively, the isomorphic classes of algebras of type
$A(1,\l,\m)$ over $D$ are in bijective correspondence with the elements of $\bar{\mathcal{M}}$.\medskip
\item\label{mohon} Let $D$ be a domain and $\Omega_3:=\{(\epsilon,\nu,\rho)\in(D^*\setminus D^\times)^3\colon\rho-\epsilon\nu\in D^\times\}$.
Consider the class of $2$-dimensional evolution algebras (over the domain $D$)
given by 
$C=\{B(1,\xi,\nu,\rho)\colon (\xi,\nu,\rho)\in\Omega_3\}$ where the product in $B(1,\xi,\nu,\rho)$ relative to a natural basis $\{e_1,e_2\}$ is 
$e_1^2=e_1+\xi e_2$, $e_2^2=\nu e_1+\rho e_2$. The isomorphism condition within
the class $C$ is $B(1,\xi,\nu,\rho)\cong B(1,\xi',\nu',\rho')$ if and only
if there is some $k\in K^\times$ such that $\xi'=\frac{1}{k}\xi$,
$\nu'=k^2\nu$ and $\rho'=k\rho$. There is an action 
$D^\times\times \Omega_3\to \Omega_3$ given by 
$k\cdot (\xi,\nu,\rho)=(\frac{\xi}{k},k^2 \nu,k \rho)$ for $k\in D^\times$.
The set of orbits of $\Omega_3$ under this action will be denoted $\Omega_3/D^\times$. Then, the isomorphic classes in $C$ are in one-to-one correspondence with the elements of the orbit set $\Omega_3/D^\times$. 
Alternatively, we can consider the curve $xz=\xi\rho$, $\rho^2y=\nu z^2$ in $Q^3$ so that the  isomorphic classes of algebras $B(1,x,y,z)$ in $C$ is in one-to-one correspondence with the points of intersection of $\Omega_3$ with the curve $$\begin{cases} xz=\xi \rho\\ \rho^2y=\nu z^2\end{cases}\qquad (\xi,\nu,\rho)\in\Omega_3.$$ Denote by $\partial_{\xi,\nu,\rho}$ such curve in the affine space $Q^3$ and $\partial_{\xi,\nu,\rho}^*:=\partial_{\xi,\nu,\rho}\cap\Omega_3$ its section with
$\Omega_3$. Then we have (as in item (I)) that $(\xi',\nu',\rho')\notin\partial_{\xi,\nu,\rho}$ implies  $\partial_{\xi,\nu,\rho}^*\cap\partial_{\xi',\nu',\rho'}^*=\emptyset$. Also $(\xi,\nu,\rho)\in\partial_{\xi,\nu,\rho}$ for any $(\xi,\nu,\rho)\in\Omega_3$. Thus, $\Omega_3$ is a disjoint union of sections $\partial_{\xi,\nu,\rho}^*$ and we can define $\bm{\partial}$ as the set whose elements are the different curves $\partial_{\xi,\nu,\rho}$. Therefore, we get a bundle 
$\Omega_3\to \bm{\partial}$ such that $(\xi,\nu,\rho)\mapsto\partial_{\xi,\nu,\rho}$ and the isomorphic
class of algebras in $C$ are in one-to-one correspondence with $\bm{\partial}$. 
So, the isomorphic classes of algebras in $C$ are indexed by the set of curves $\bm{\partial}$ or if we prefer, a moduli set for $C$ is $\bm{\partial}$: each algebra in $C$ is completely determined by a curve $\partial_{\xi,\nu,\rho}$.
\item\label{elragep}
Again let $D$ denote a domain with field of fractions $Q$ and define $\Sigma_3:=\{(\m,\l,\o)\in (D^*\setminus D^\times)^3\colon \m \o -\l\in D^\times\}$.
Now let $C$ denote the class of algebras $B(\m,1,\l,\o)$ whose multiplication in
a natural basis is $e_1^2=\m e_1+e_2$, $e_2^2=\l e_1+\o e_2$.
We have an action $D^\times\times\Sigma_3\to\Sigma_3$ given by 
    $k(\m,\l,\o)=(k\m,k^3\l,k^2\o)$. The isomorphic classes of algebras of this kind are in one-to-one correspondence with the orbit set $\Sigma_3/D^\times$. Let $\sigma_{\m,\l,\o}$ be the curve of $Q^3$ given by
    \begin{equation}\label{narram}
    \begin{cases}\mu^3y=x^3\lambda\\\mu^2z=x^2\o
    \end{cases} \qquad (\m,\l,\o)\in\Sigma_3.
    \end{equation}
    Denote as in previous cases $\s_{\m,\l,\o}^*=\s_{\m,\l,\o}\cap\Sigma_3$.
    Then, $(\m,\l,\o)\in \s_{\m,\l,\o}$ for any $(\m,\l,\o)\in \Sigma_3$. As before, the different $\s_{\m,\l,\o}^*$'s are pairwise disjoint and their disjoint union is $\Sigma_3$ so, we get a bundle $\Sigma_3\to\bm{\s}$ where $\bm{\s}$ is the set whose elements are the $\s_{\m,\l,\o}$'s.
    Thus, $\Sigma_3/D^\times$ is is one-to-one correspondence with the set $\bm{\s}$ and this is a moduli for the class $C$.
\item\label{mero} 
For a  domain $D$, with field of fractions $Q$, define the group $\m_3(D):=\{k\in D^\times \colon k^3=1\}$ and
    $\Omega_2:=\{(\xi,\rho)\in (D^*\setminus D^\times)^2\colon \xi \rho-1\in D^\times\}$.
    We have an action $\m_3(D)\times\Omega_2\to\Omega_2$ given by 
    $k\cdot(\xi,\rho)=(k\xi,k^{-1}\rho)$. To analyze the orbit set $\Omega_2/\m_3(D)$ define for any $(\xi,\rho)\in\Omega_2$ the curve
$h_{\xi,\rho}$ of $Q^2$ given by $xy=\xi \rho$. Define also 
$h_{\xi,\rho}^*=h_{\xi,\rho}\cap\Omega_2$. We have $(\xi,\rho)\in h_{\xi,\rho}$ and if $(\xi',\rho')\notin h_{\xi,\rho}$ then 
$h_{\xi,\rho}^*\cap h_{\xi',\rho'}^*=\emptyset$ so, $\Omega_2$ is
the disjoint union of all $h_{\xi,\rho}^*$. Therefore, 
$\Omega_2/\m_3(D)$ is in one-to-one correspondence with the set $\bm{h}$ whose elements are the curves $h_{\xi,\rho}$.
\item\label{eumsut}
With $D$ and $Q$ as in previous cases, consider the set $S$ of all triples 
    $$(\m,\l,\o)\in (D^*\setminus D^\times)\times (D^\times\setminus(D^\times)^{[3]})\times (D^*\setminus D^\times),$$ such that $\ \m \o-\l\in D^\times$ and the action $D^\times\times S\to S$ given by $k(\m,\l,\o)=(k\m,k^3\l,k^2\o)$. Thus, the orbit set $S/D^\times$ is in one-to-one correspondence with the set $\bm{\s'}$ whose elements are the sections $\tilde\sigma_{\m,\l,\o}:=\sigma_{\m,\l,\o}\cap S$ defined for $(\m,\l,\o)\in S$ in the affine space $Q^3$ by equations \eqref{narram} by replacing $(\m,\l,\o)\in\Sigma_3$ with $(\m,\l,\o)\in S$.
\item\label{derniere}
We will denote $\Omega_4$ as the set of all 
$(\a,\b,\g,\d)\in (D^*\setminus D^\times)^4$ such that
$\a \d- \g \b\in D^\times$.
We have an action
$$(D^\times\times D^\times)\times\Omega_4\to\Omega_4$$
given by $(k_1,k_2)\cdot(\a,\b,\g,\d):=(k_1\a, \frac{k_1^2}{k_2}\b,\frac{k_2^2}{k_1}\g,k_2\d)$. The  orbit set $\frac{\Omega_4}{(D^\times\times D^\times)}$ can be described defining (for every $(\a,\b,\g,\d)\in\Omega_4$) the surface $\o_{\a,\b,\g,\d}$ of $Q^4$ given by
$$\begin{cases}
yt\a^2=\b \d x^2\\ z x\d^2=\a \g t^2.
\end{cases}$$
Therefore, $\bar\o_{\a,\b,\g,\d}:=\o_{\a,\b,\g,\d}\cap\Omega_4$ so that $\Omega_4$ is a disjoint union of $\bar\o_{\a,\b,\g,\d}$'s and we have a bijection
$\dfrac{\Omega_4}{D^\times\times D^\times}\cong\bm{\o}$ where the elements of $\bm{\o}$
are the sections $\bar\o_{\a,\b,\g,\d}$.
Furthermore, this orbit set is in one-to-one correspondence with
the $Q$-points of the surface
$$\begin{cases}
yt\a^2=\b \d x^2\\ z x\d^2=\a \g t^2
\end{cases}$$
which are in $\Omega_4$.
\end{enumerate}

\section{The perfect case}\label{perfectcase}

In the previous section we have analyzed some of the different moduli sets that we will use  now  in our  classification task. Consider a $2$-dimensional evolution algebra $\E$ over the domain $D$ with a natural basis $\{e_1,e_2\}$ and assume that $\E$ is perfect.
We analyse several cases.

\subsection{Case $e_1^2=\a e_1$,  $e_2^2=\b e_2$}
Then $\a,\b\in D^\times$. We can define 
$f_1=\a^{-1}e_1$ and $f_2=\b^{-1}e_2$. So we have 
$f_1^2=\a^{-2}\a e_1=f_1$ and 
$f_2^2=\b^{-2}\b e_2=f_2$. All the algebras in this case are isomorphic to $D\times D$ with componentwise operations. Then, up to isomorphism, there is only one algebra of this type.

\subsection{Case $e_1^2=\a e_2$, $e_2^2=\b e_1$}
Again $\a,\b\in D^\times$. We can define 
$f_1=e_1$ and $f_2=\a e_2$. So, we have 
$f_1^2=\a e_2=f_2$ and $f_2^2=\a^2\b f_1$. Thus, we have
a one-parametric family of algebras given by the multiplication table $e_1^2=e_2$ and $e_2^2=\a e_1$ where $\a\in D^\times$. Denote the above algebra by $A_{2,\a}$. We analyse the isomorphism question
$A_{2,\a}\cong A_{2,\b}$. If the isomorphism exists we have bases $\{e_1,e_2\}$ and $\{u_1,u_2\}$ of $A_{2,\a}$
such that $e_1^2=e_2$, $e_2^2=\a e_1$ and 
$u_1^2=u_2$, $u_2^2=\b u_1$. Then, we have two possibilities:
\begin{enumerate}
    \item $u_i=k_i e_i$, $k_i \in D^{\times}$ for $i=1,2$. This gives 
    $\b\a^{-1}=k_1^3$ and $k_2=k_1^2$. So $\b\a^{-1}\in (D^\times)^3$.
    \item $u_1=k e_2$ and $u_2=h e_1$ for some $k,h\in D^\times$. In this case $k^3=\b\a^{-2}$, $h=k^2\a$ and so
    $\b\a^{-2}\in (D^\times)^3.$
\end{enumerate}
Therefore, we have $A_{2,\a}\cong A_{2,\b}$ if and only if 
$\b\a^{-1}\in(D^\times)^3$ or 
$\b\a^{-2}\in(D^\times)^3$. But, I ask the reader to check the implication $\Leftarrow$. A moduli set for the class of algebras $A_{2,\a}$ is the group
$\lim_{\to 2}D^\times/(D^\times)^{[3]}$.

\subsection{Case $e_1^2=\a e_1$, $e_2^2=\b e_1+\d e_2$, $\b\ne 0$}
In this case $\a,\d\in D^\times$. We define 
$f_1=\a^{-1}e_1$ and $f_2=\d^{-1}e_2$. Then 
$f_1^2=f_1$ and $f_2^2=\l f_1+f_2$ for a certain $\l\in D^*$. We get a one-parameter family of algebras $A_{3,\l}$ with $\l\ne 0$ and product $e_1^2=e_1$, $e_2^2=\l e_1+e_2$. We investigate the isomorphism $A_{3,\l}\cong A_{3,\m}$. As before, we have two natural basis  $\{e_1,e_2\}$ and $\{u_1,u_2\}$ such that 
$e_1^2=e_1$, $e_2^2=\l e_1+e_2$ and 
$u_1^2=u_1$, $u_2^2=\m u_1+u_2$. Then:
\begin{enumerate}
    \item If $u_i=k_i e_i$ for $i=1,2$. After some computations, we get 
    $k_2=k_1=1$, thus $\l=\m$.
    \item If $u_1=k e_2$ and $u_2=h e_1$ we get an
    inconsistent system of equations.
\end{enumerate}
Thus, $A_{3,\l}\cong A_{3,\m}$ if and only if $\l=\m$.
The moduli set for the class $A_{3,\l}$ is $D^*$ which is a monoid. 
\subsection{Case $e_1^2=\a e_2$, $e_2^2=\b e_1+\d e_2$, $\d\ne 0$}
In this case, $\a,\b\in D^\times$. We have a family
$A(\a,\b,\d)$ of algebras with $\a,\b\in D^\times$ and 
$\d \in D^*$ (and the multiplication above). It is straightforward to prove that there is no possible isomorphism 
$A(\a,\b,\d)\cong A(\a',\b',\d')$ when $\d\in D^\times$ but $\d'\not\in D^\times$. So, let us investigate the algebras $A(\a,\b,\d)$ with $\d\in D^\times$. We will prove that in this case $A(\a,\b,\d)\cong A(\l,1,1)$ for a suitable $\l\in D^\times$. Indeed, defining $f_1=\b\d^{-2}e_1$
and $f_2=\d^{-1}e_2$ one gets $f_1^2=\l f_2$ and 
$f_2^2=f_1+f_2$ where $\l=\a\b^2\d^{-3}\in D^\times$.
Furthermore, it is easy to check that $A(\l,1,1)\cong A(\mu,1,1)$ if and only if $\l=\mu$. Let us investigate now the other class of algebras: $A(\a,\b,\d)$ with $\d\not\in D^\times$. By making the change of basis
$f_1=k_1e_1$, $f_2=\a k_1^2 e_2$ we get 
$f_1^2=k_1^2\a e_2=f_2$ and 
$f_2^2=\a^2\b k_1^3 f_1+\d\a k_1^2 f_2$ and $\a^2\b k_1^3\in D^\times$ while 
$\d\a k_1^2\not\in D^\times$. Thus, any $A(\a,\b,\d)$ with $\d$ not invertible is
isomorphic to $A(1,\l,\mu)$ where $\l\in D^\times$ and $0\neq \mu\not\in D^\times$. In addition, $A(1,\l,\mu)\cong A(1,k^3\l,k^2\mu)$ for any $k\in D^\times$ 
($\mu\notin D^\times$). Moreover, $A(1,\l,\mu)\cong A(1,\l',\mu')$ (where $\mu,\mu'\not\in D^\times$) if and only if there is a $k\in D^\times$ such that 
$\l'=k^3\l$ and $\mu'=k^2\mu$.
Summarizing: the algebras in this case fall into two mutually non-isomorphic classes: those of the form $A(\lambda,1,1)$ with $\l\in D^\times$ and those of the form $A(1,\l,\mu)$ with $\l\in D^\times$, $0 \neq \mu\notin D^\times$. Also,
$$\begin{matrix}
A(\l,1,1)\cong A(\mu,1,1) & \hbox{ iff } &
 \l=\mu\cr 
 A(1,\l,\mu)\cong A(1,\l',\mu')  & \hbox{ iff } & \exists k\in D^\times\colon \l'=k^3\l, \mu'=k^2\mu.\end{matrix}$$
 The algebras of the form $A(1,\l,\mu)$ only exist over domains which are not fields.
 The isomorphic classes of algebras of type $A(1,\l,\m)$ are in one-to-one correspondence with the monoid $\bar{\mathcal{M}}$ defined in equation \eqref{sapr1}.

\subsection{Case $e_1^2=\a e_1+\b e_2$, $e_2^2=\g e_1+\d e_2$, $\a,\b,\g,\d\ne 0$}
In this case $\a\d-\b\g\in D^\times$. Denote by 
$B(\a,\b,\g,\d)$ the two-dimensional evolution algebra with natural basis $\{e_1,e_2\}$ and multiplication $e_1^2=\a e_1+\b e_2$, $e_2^2=\g e_1+\d e_2$ being $\a\d-\b\g\in D^\times$ and
$\a,\b,\g,\d\ne 0$. The change 
$f_i=k_i e_i$ ($i=1,2$) and $k_i\in D^\times$ gives 
$$\begin{cases}f_1^2=k_1\a f_1+\frac{k_1^2}{k_2}\b f_2\cr 
f_2^2=\frac{k_2^2}{k_1}\g f_1+k_2\d f_2.\end{cases}$$
So,
\begin{equation}\label{ecuacion1}
    B(\a,\b,\g,\d)\cong B(k_1\a,\frac{k_1^2}{k_2}\b,\frac{k_2^2}{k_1}\g,k_2\d)
\end{equation}
for any $k_i\in D^\times$. On the other hand, the change $f_1=ke_2$, $f_2=h e_1$ with $k,h\in D^\times$ produces
$$\begin{cases}f_1^2=k\d f_1+\frac{k^2\g}{h}f_2\cr 
f_2^2=\frac{h^2\b}{k} f_1+h\a f_2\end{cases}$$
which allows to conclude that 
 $B(\a,\b,\g,\d)\cong B(k\d,\frac{k^2\g}{h},\frac{h^2\b}{k},h\a)$ for any $k,h\in D^\times$. Any isomorphism between algebras of the type $B(\a,\b,\g,\d)$ is of one of the previous forms.
 We now distinguish several mutually non-isomorphic classes:

\subsubsection{Both $\a$ and $\d$ are invertible}
 We have the following isomorphism $B(\a,\b,\g,\d)\cong B(1,\frac{\d}{\a^2}\b,\frac{\a}{\d^2}\g,1)$. Thus, our algebra is of type $B(1,\l,\m,1)$ with $1-\l \m \in D^\times$ and $\l,\m\ne 0$. Moreover, it is easy to see that $B(1,\l,\m,1)\cong B(1,\l',\m',1)$ if and only if $(\l,\m)=(\l',\m')$ or $(\l,\m)=(\m',\l')$. The moduli set is the orbit set $X/\Z_2$ where $X$ is the set of all $(\l,\m)\in D^*\times D^*$ with $1-\l \m \in D^\times$ and the action $\Z_2\times X\to X$ is 
$0\cdot (\l,\m)=(\l,\m)$, $1\cdot(\l,\m)=(\m,\l)$.

\subsubsection{Only one of $\a$ and $\d$ is invertible}
We may assume without loss of generality that $\a\in D^\times$ but $\d\notin D^\times$.
Thus, our algebra is isomorphic to some 
$B(1,\xi,\nu,\rho)$ with $\rho-\xi \nu \in D^\times$, $\xi,\nu,\rho\ne 0$ and $\rho\notin D^\times$. We focus then in the class of algebras $${\mathcal C}:=\{B(1,\xi,\nu,\rho)\colon\ 
\rho-\xi \nu \in D^\times, \xi,\nu, \rho \ne 0, \rho \notin D^\times\}.$$
But then, ${\mathcal C}={\mathcal C}_1\sqcup {\mathcal C}_2\sqcup {\mathcal C}_3$ (a disjoint union) where 
\begin{eqnarray*}
{\mathcal C}_1&:= &\{B(1,\xi,\nu,\rho)\in\mathcal{C}\colon \xi,\nu\in D^\times, \rho \notin D^\times\},\cr
{\mathcal C}_2&:= & \{B(1,\xi,\nu,\rho)\in\mathcal{C}\colon \xi\in D^\times, \nu, \rho \notin D^\times\},\cr
{\mathcal C}_3&:= & \{B(1,\xi,\nu,\rho)\in\mathcal{C}\colon \nu\in D^\times, \xi, \rho \notin D^\times\},\cr
{\mathcal C}_4&:= &\{B(1,\xi,\nu,\rho)\in\mathcal{C}\colon \xi,\nu, \rho\notin D^\times\}.\end{eqnarray*}

The algebras of the different classes are not isomorphic and any algebra in $\mathcal C$ is isomorphic to some of ${\mathcal C}_i$, for $i\in\{1,2,3,4\}$.
Thus, we have to investigate the isomorphism question within each class ${\mathcal C}_i$. 
\begin{enumerate}
\item[${\mathcal C}_1$:] If $k\in D^\times$ we have 
$B(1,\xi,\nu,\rho)\cong B(1,\frac{\xi}{k},k^2 \nu, k\rho)$ using \eqref{ecuacion1} and for $k=\xi$ we have 
$B(1,\xi,\nu,\rho)\cong B(1,1,\xi^2 \nu, \xi \rho)$. Then, any algebra in ${\mathcal C}_1$ is isomorphic so $B(1,1,\lambda,\mu)$ with $\lambda\in D^\times$, $\mu\notin D^\times$, $\mu\ne 0$ and $\lambda-\mu\in D^\times$. Moreover, two such algebras $B(1,1,\lambda,\mu)$ and 
$B(1,1,\lambda',\mu')$ are isomorphic if and only if $(\lambda,\mu)=(\lambda',\mu')$. For instance, over the integers there are only two non-isomorphic classes in ${\mathcal C}_1$: the one of $B(1,1,1,2)$ and that of $B(1,1,-1,-2)$.

\item[${\mathcal C}_2$:] Any algebra in this class is isomorphic to some $B(1,1,\l,\m)$ with $\l,\m\notin D^\times$ but $\l-\m\in D^\times$ and $\l,\m\ne 0$.
For any such algebras we have $B(1,1,\l,\m)\cong B(1,1,\l',\m')$ if and only if $(\l,\m)=(\l',\m')$.

\item[${\mathcal C}_3$:] In this case, any algebra in this class is isomorphic to some $B(1,\l,1,\m)$ with $\l,\m\notin D^\times$, but $\m-\l\in D^\times$ and $\l,\m\ne 0$.
For any such algebras we have $B(1,\l,1,\m)\cong B(1,\l',1,\m')$ if and only if $(\l,\m)=(\l',\m')$ or $(\l,\m)=(-\l',-\m')$.

\item[${\mathcal C}_4$:] We have $B(1,\xi,\nu,\rho)\cong B(1,\frac{\xi}{k},k^2 \nu, k\rho)$ for $k\in D^\times$ using \eqref{ecuacion1}. Denote by $D^*\setminus D^\times$ the set of nonzero and noninvertible elements of $D$ and recall that
$$\Omega_3=\{(\xi,\nu,\rho)\in (D^*\setminus D^\times)^3\colon \rho-\xi \nu\in D^\times\}.$$ 
The isomorphism condition in this class of algebras is given in \ref{tnerap} (II) and
the moduli set is that of the curves $\partial_{\xi,\nu,\rho}$.

For instance, if $D={\mathbb Z}$, an element in this class of algebras is $B(1,3,2,5)$ and its isomorphic class consists on itself and  $B(1,-3,2,-5)$.

\end{enumerate}

\subsubsection{Neither $\a$ nor $\d$ are invertible}
Again we distinguish two mutually non-isomorphic cases.
\begin{enumerate}
    \item\label{griego} $\beta$ or $\g$ is invertible. We may assume without loss of generality that $\b\in D^\times$. Then $B(\a,\b,\g,\d)\cong B(k\a,1,\b^2k^3\g,\b\d k^2)$  for any $k\in D^\times$ using \eqref{ecuacion1}. The algebras in this class are therefore of the form $B(\m,1,\l,\o)$ with $\m,\o\notin D^\times$ and 
    $\m,\l,\o\ne 0$, $\m \o-\l\in D^\times$. 
    \begin{enumerate}
    \item Assume $\l \notin D^\times$. Since $B(\m,1,\l,\o)\cong B(k\m,1,k^3\l,k^2\o)$ for any $k\in D^\times$, 
    we may define $\Sigma_3:=\{(\m,\l,\o)\in (D^*\setminus D^\times)^3\colon \m \o -\l\in D^\times\}$ and
    we have an action $D^\times\times\Sigma_3\to\Sigma_3$ given by 
    $k(\m,\l,\o)=(k\m,k^3\l,k^2\o)$. The isomorphic classes of algebras of this kind are in one-to-one correspondence with the orbit set $\Sigma_3/D^\times$. Alternatively, $\Sigma_3/D^\times$ is is one-to-one correspondence with the set of curves $\bm{\sigma}$ defined in \ref{tnerap} (III).
    \item\label{musa} If $\l\in (D^\times)^{[3]}$, then $\l=\epsilon^3$ for some invertible $\epsilon$. Since $B(\m,1,\l,\o)\cong B(k\m,1,k^3\l,k^2\o)$ (for $k\in D^\times$), we may take $k=\epsilon^{-1}$ and then 
    $k^3 \l=1$ so that $B(\m,1,\l,\o)\cong B(k\m,1,1,k^2\o)$. Thus, the algebras in this case are all of the form $B(\xi,1,1,\rho)$ with $\xi,\rho\ne 0, \xi \rho-1\in D^\times$ and $\xi,\rho\notin D^\times$. We have 
    $B(\xi,1,1,\rho)\cong B(k\xi,1,1,k^{-1}\rho)$ whenever
    $k\in D^\times$ satisfies $k^3=1$. So, denoting $\m_3(D):=\{k\in D^\times \colon k^3=1\}$, and
    $\Omega_2:=\{(\xi,\rho)\in (D^*\setminus D^\times)^2\colon \xi \rho-1\in D^\times\}$,
    we have the action defined in \ref{tnerap} (IV) given by $\m_3(D)\times\Omega_2\to\Omega_2$ with 
    $k\cdot(\xi,\rho)=(k\xi,k^{-1}\rho)$. Consequently, the isomorphic classes of algebras in this case are in one to one correspondence with the orbit set $\Omega_2/\m_3(D)$ which can be identified with the set $\bm{h}$ of \lq\lq hyperbolae\rq\rq\ $h_{a,b}$ (see \ref{tnerap} (IV)).    
    \item\label{orom} In case $\l\in D^\times$ but not necessarily $\l\in (D^\times)^{[3]}$, 
    since $B(\m,1,\l,\o)\cong B(k\m,1,k^3\l,k^2\o)$ for any $k\in D^\times$,
    the only thing we can do is to consider the set $S$ of all triples 
    $(\m,\l,\o)\in (D^*\setminus D^\times)\times (D^\times\setminus(D^\times)^{[3]})\times (D^*\setminus D^\times),$ such that $\ \m \o-\l\in D^\times$ and the action $D^\times\times S\to S$ given by $k(\m,\l,\o)=(k\m,k^3\l,k^2\o)$. Then, the isomorphic classes of this type are in one-to-one correspondence with the elements of $S/D^\times$.  This can be identified with the set of sections $\bm{\sigma}'$ defined in \ref{tnerap} (V).
\end{enumerate}
Note that the case presented in (\ref{musa}), that is, when $\l \in (D^\times)^{[3]}$, is in fact a subcase of (\ref{orom}). We have specified it because if $\l$ is a cube, one more $1$ can be got in the structure matrix of the algebra.
\item $\b,\g\notin D^\times$. So, we have the algebras 
$B(\a,\b,\g,\d)$ where the four scalars are nonzero and noninvertible but $\a\d-\b\g\in D^\times$.
We will denote $\Omega_4$ as the set of all 
$(\a,\b,\g,\d)\in (D^*\setminus D^\times)^4$ such that
$\a \d- \g \b\in D^\times$.
We have an isomorphism $B(\a,\b,\g,\d)\cong B(k_1\a,\frac{k_1^2}{k_2}\b,\frac{k_2^2}{k_1}\g,k_2\d)$
for any $k_1,k_2\in D^\times$. Thus we have an action
$$(D^\times\times D^\times)\times\Omega_4\to\Omega_4$$
given by $(k_1,k_2)\cdot(\a,\b,\g,\d):=(k_1\a, \frac{k_1^2}{k_2}\b,\frac{k_2^2}{k_1}\g,k_2\d)$. The isomorphic classes of algebra of this type are in one-to-one correspondence with the elements of the orbit set $\Omega_4/(D^\times\times D^\times)$. 
Furthermore, this orbit set is in one-to-one correspondence with the set of sections $\bm{\o}$ described in \ref{tnerap} (VI).
For instance, if $D=\mathbb{Z}$ the isomorphic class of $B(2,3,3,5)$ has the following four elements:
$B(2,3,3,5)$, $B(-2,3,-3,5)$, $B(2,-3,3,-5)$ and $B(-2,-3,-3,-5)$. 
Now, if we consider $D={\mathbb Z}[x,x^{-1}]$ then $D^\times=\{\pm x^n\colon n\in{\mathbb Z}\}$ and we have 
$$B(7x^3+4x^2,4x,5x^5+3x^4,3x^3)\cong B(7x+4,4,5x+3,3)$$ being 
the elements $(7x^3+4x^2,4x,5x^5+3x^4,3x^3)$ and $(7x+4,4,5x+3,3)$ in the same orbit of $\Omega_4$ under the action of $D^\times\times D^\times$.
Indeed, $(x^2,x^3)\cdot (7x+4,4,5x+3,3)=(7x^3+4x^2,4x,5x^5+3x^4,3x^3).$
\end{enumerate}

Now, we will collect all these results in the following theorem.

\begin{theorem}\label{teoremon}
Let $\E$ be a perfect two-dimensional evolution algebra over a domain $D$. Then, we have one and only one of the following possibilities:
\begin{enumerate}
    \item $\E\cong A_1=D\times D$ with product $(x,y)(u,v)=(xu,yv)$. The graph associated to the evolution algebra is of the form
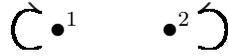
\[
\xymatrix{\bullet^{1} \ar@(dl,ul) \ar@(dl,ul)@<0.2mm>   & \bullet^{2} \ar@(dr,ur) \ar@(dr,ur)@<0.2mm>}\]

    \item $\E\cong A_{2,\a}$ where $A_{2,\a}=D\times D$ with product $(x,y)(u,v)=(\a y v,x u)$ and $\a\in D^\times$. This algebras are classified by the moduli $\lim_{\to 2}D^\times/(D^\times)^{[3]}$, that is;  $A_{2,\a}\cong A_{2,\b}$ if and only if $\b\a^{-1}\in (D^\times)^3$ or $\b\a^{-2}\in (D^\times)^3$. The graph is
$$
\xymatrix{\bullet^{1}  \ar@/^.5pc/[r] \ar@<0.2mm>@/^.5pc/[r] & \bullet^{2} \ar@/^.5pc/[l] \ar@<0.2mm>@/^.5pc/[l]  }$$

    \item $\E\cong A_{3,\a}=D\times D$ with product $(x,y)(u,v)=(xu+\a yv,yv)$ and $\a\ne 0$. This algebras are classified by the moduli $D^*$, that is $A_{3,\a}\cong A_{3,\b}$ if and only if $\a=\b$. The graph  is one of

$$
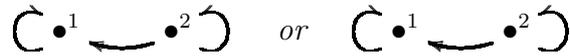
\xymatrix{\bullet^{1} \ar@(dl,ul) \ar@<0.2mm>@(dl,ul) & \bullet^{2} \ar@(dr,ur) \ar@<0.2mm>@(dr,ur) \ar@[blue]@/^.5pc/[l]
\ar@[blue]@<0.2mm>@/^.5pc/[l]  }
\qquad\hbox{ or }
\qquad
\xymatrix{\bullet^{1} \ar@(dl,ul) \ar@<0.2mm>@(dl,ul) & \bullet^{2} \ar@(dr,ur) \ar@<0.2mm>@(dr,ur) \ar@/^.5pc/[l]  \ar@<0.2mm>@/^.5pc/[l]  \ar@/^.5pc/[l]}
$$
\smallskip

\noindent where the blue arrow stands for the fact that $\a$ is not invertible but it is nonzero. 
    \item Denote by $A(\a,\b,\d)=D\times D$ with multiplication 
    $(x,y)(u,v)=(\b yv,\d yv+\a xu)$. Then,
\begin{enumerate}    
    \item $\E\cong A(\l,1,1)=D\times D$ with $\l\in D^\times$. The corresponding graph  is 
$$
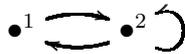
\xymatrix{\bullet^{1} \ar@/^.5pc/[r] \ar@<0.2mm>@/^.5pc/[r]     & \bullet^{2} \ar@(dr,ur) \ar@<0.2mm>@(dr,ur) \ar@/^.5pc/[l] \ar@<0.2mm>@/^.5pc/[l] }$$
    \item $\E\cong A(1,\l,\m)$ with $\l\in D^\times$ and $0\ne\m\notin D^\times$. The graph is
    
$$
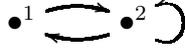
\xymatrix{\bullet^{1} \ar@/^.5pc/[r] \ar@<0.2mm>@/^.5pc/[r] & \bullet^{2} \ar@[blue]@(dr,ur) \ar@<0.2mm>@[blue]@(dr,ur)  \ar@/^.5pc/[l] \ar@<0.2mm>@/^.5pc/[l]   }
$$
\end{enumerate}
An algebra in one of the cases is not isomorphism to any algebra in the other case. The isomorphism condition in each case are: $A(\l,1,1)\cong A(\m,1,1)$ if and only if $\l=\m$, that is, the moduli is $D^\times$; and 
$A(1,\l,\mu)\cong A(1,\l',\m')$ if and only if $\exists k\in D^\times$, $\l'=k^3\l$, $\m'=k^2 \m$. The moduli set is $\bar{\mathcal{M}}$.
\item Denote by $B(\a,\b,\g,\d)=D\times D$ with multiplication $(x,y)(u,v)=(\a xu+\g yv,\b xu+\d yv)$. Then: 
\begin{enumerate}
\item[(5.I)] Either $\a$ or $\d$ is invertible. We have the following mutually excluding cases:
\begin{enumerate} 
    \item  $\E\cong B(1,\l,\m,1)$, $\l,\m\ne 0$, $1-\l\m\in D^\times$. The graph is one of
       $$\qquad \qquad 
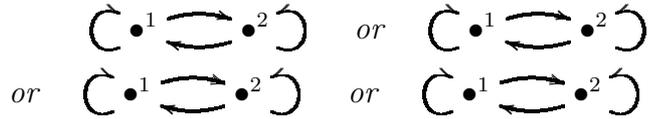
\xymatrix{\bullet^{1} \ar@(dl,ul) \ar@<0.2mm>@(dl,ul) \ar@/^.5pc/[r]  \ar@<0.2mm>@/^.5pc/[r]  & \bullet^{2} \ar@(dr,ur) \ar@<0.2mm>@(dr,ur) \ar@/^.5pc/[l]  \ar@<0.2mm>@/^.5pc/[l]  }
    \qquad\hbox{ or }
\qquad
    \xymatrix{\bullet^{1} \ar@(dl,ul) \ar@<0.2mm>@(dl,ul) \ar@[blue]@/^.5pc/[r]  \ar@[blue]@<0.2mm>@/^.5pc/[r]  & \bullet^{2} \ar@(dr,ur) \ar@<0.2mm>@(dr,ur) \ar@/^.5pc/[l]  \ar@<0.2mm>@/^.5pc/[l] }
     $$
   $$ \hbox{ or } \qquad \xymatrix{\bullet^{1} \ar@(dl,ul) \ar@<0.2mm>@(dl,ul) \ar@/^.5pc/[r]  \ar@<0.2mm>@/^.5pc/[r]  & \bullet^{2} \ar@(dr,ur) \ar@<0.2mm>@(dr,ur) \ar@[blue]@/^.5pc/[l]  \ar@[blue]@<0.2mm>@/^.5pc/[l] }
    \qquad\hbox{ or }
\qquad
\xymatrix{\bullet^{1} \ar@(dl,ul) \ar@<0.2mm>@(dl,ul) \ar@[blue]@/^.5pc/[r]  \ar@[blue]@<0.2mm>@/^.5pc/[r]  & \bullet^{2} \ar@(dr,ur) \ar@<0.2mm>@(dr,ur) \ar@[blue]@/^.5pc/[l]  \ar@[blue]@<0.2mm>@/^.5pc/[l]  } $$
        \item $\E\cong B(1,1,\l,\m)$ with $\l\in D^\times$, $0\ne \m\notin D^\times$, $\l-\m\in D^\times$. The corresponding graph is     $$\qquad \qquad 
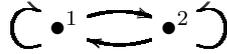
\xymatrix{\bullet^{1} \ar@(dl,ul) \ar@<0.2mm>@(dl,ul) \ar@/^.5pc/[r]  \ar@<0.2mm>@/^.5pc/[r]  & \bullet^{2} \ar@[blue]@(dr,ur) \ar@[blue]@<0.2mm>@(dr,ur) \ar@/^.5pc/[l]  \ar@<0.2mm>@/^.5pc/[l]}$$ \label{prueba1}
    \item $\E\cong B(1,1,\l,\m)$ with $\l,\m\notin D^\times$ but $\l-\m\in D^\times$, $\l,\m\ne 0$. The graph is
    
        $$\qquad \qquad
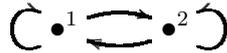
\xymatrix{\bullet^{1} \ar@(dl,ul) \ar@<0.2mm>@(dl,ul) \ar@/^.5pc/[r]  \ar@<0.2mm>@/^.5pc/[r]  & \bullet^{2} \ar@[blue]@(dr,ur) \ar@[blue]@<0.2mm>@(dr,ur) \ar@[blue]@/^.5pc/[l]  \ar@[blue]@<0.2mm>@/^.5pc/[l] }$$ \label{prueba2}
        
         \item $\E\cong B(1,\l,1,\m)$ with $\l,\m\notin D^\times$ but $\m-\l\in D^\times$, $\l,\m\ne 0$. The associated graph is
         
         $$ \qquad \qquad
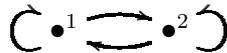
\xymatrix{\bullet^{1} \ar@(dl,ul) \ar@<0.2mm>@(dl,ul) \ar@[blue]@/^.5pc/[r]  \ar@[blue]@<0.2mm>@/^.5pc/[r]  & \bullet^{2} \ar@[blue]@(dr,ur) \ar@[blue]@<0.2mm>@(dr,ur) \ar@/^.5pc/[l]  \ar@<0.2mm>@/^.5pc/[l] }$$
        \item $\E\cong B(1,\l,\m,\g)$, $\l,\m,\g\ne 0$, $\g-\l\m\in D^\times$, $\g,\l,\m\notin D^\times$. The graph is
        $$\qquad \qquad
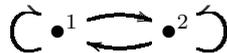
\xymatrix{\bullet^{1} \ar@(dl,ul) \ar@<0.2mm>@(dl,ul) \ar@[blue]@/^.5pc/[r]  \ar@[blue]@<0.2mm>@/^.5pc/[r]  & \bullet^{2} \ar@[blue]@(dr,ur) \ar@[blue]@<0.2mm>@(dr,ur) \ar@[blue]@/^.5pc/[l]  \ar@[blue]@<0.2mm>@/^.5pc/[l]  }$$
\end{enumerate}
An algebra in one of the cases is not isomorphism to any algebra in other case.
In case (i), $B(1,\l,\m,1)\cong B(1,\l',\m',1)$ if and only if $(\l,\m)=(\l',\m')$ or $(\l,\m)=(\m',\l')$. 
In cases (ii) and (iii),   $B(1,1,\l,\m)\cong B(1,1,\l',\m')$ if and only if $(\l,\m)=(\l',\m')$. In case (iv) $B(1,\l,1,\m)\cong B(1,\l',1,\m')$ if and only if $(\l,\m)=(\l',\m')$ or $(\l,\m)=(-\l',-\m')$. In case (v) $B(1,\l,\m,\g)\cong B(1,\l',\m',\g')$ if and only if
$\exists k\in D^\times$, $\lambda'=k^{-1}\l$,
$\m'=k^2\m$ and $\g'=k\g$. The moduli set is $\bm{\partial}$.
\item[(5.II)] Neither $\a$ nor $\d$ is invertible, but $\b$ or $\d$ is invertible. We have the possibilities:
\begin{enumerate}
    \item $\E\cong B(\l,1,\m,\g)$, $\l,\g,\mu\notin D^\times$, $\l,\m,\g\ne 0$, $\g\l-\m\in D^\times$. The graph is
    
    $$\qquad \qquad
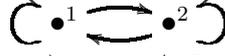
\xymatrix{\bullet^{1} \ar@[blue]@(dl,ul) \ar@[blue]@<0.2mm>@(dl,ul) \ar@/^.5pc/[r]  \ar@<0.2mm>@/^.5pc/[r]  & \bullet^{2} \ar@[blue]@(dr,ur) \ar@[blue]@<0.2mm>@(dr,ur) \ar@[blue]@/^.5pc/[l]  \ar@[blue]@<0.2mm>@/^.5pc/[l]  }$$

    \item $\E\cong B(\l,1,\m,\g)$, $\l,\g\notin D^\times$, $\mu\in D^\times$ 
    $\l,\m,\g\ne 0$, $\g\l-\m\in D^\times$. The corresponding graph is
    
     $$\qquad \qquad
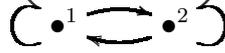
\xymatrix{\bullet^{1} \ar@[blue]@(dl,ul) \ar@[blue]@<0.2mm>@(dl,ul) \ar@/^.5pc/[r]  \ar@<0.2mm>@/^.5pc/[r]  & \bullet^{2} \ar@[blue]@(dr,ur) \ar@[blue]@<0.2mm>@(dr,ur) \ar@/^.5pc/[l]  \ar@<0.2mm>@/^.5pc/[l]  }$$

In this case if $\mu\in(D^\times)^{[3]}$, we have
$\E\cong B(\l,1,1,\g)$.

\end{enumerate}
Again the cases are mutually excluding. In case (i) we have $B(\l,1,\m,\g)\cong B(\l',1,\m',\g')$ if and only if $\exists k\in D^\times$, $\l'=k\l$, $\m'=k^3\m$, $\g'=k^2\g$.  The moduli set is $\bm{\sigma}$. 
In case (ii), when $\mu\notin(D^\times)^{[3]}$, we have $B(\l,1,\m,\g)\cong B(\l',1,\m',\g')$ if and only if 
$\exists k\in D^\times$, $\l'=k\l$, $\m'=k^3\m$, $\g'=k^2\g$. The moduli set is $\bm{\sigma}'$. In this case, when $\mu\in(D^\times)^{[3]}$, we get $B(\l,1,1,\g)\cong B(\l',1,1,\g')$ if and only if $\exists k \in D^\times$, $\l'=k \l$ and $\m'=k^{-1} \m$. The moduli set is $\bm{h}$.
\item[(5.III)] The elements $\a$, $\d$, $\b$ and $\d$ are not invertible. Then $\E\cong B(\a,\b,\g,\d)$ and the graph is

 $$\qquad \qquad
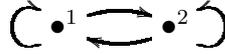
\xymatrix{\bullet^{1} \ar@[blue]@(dl,ul) \ar@[blue]@<0.2mm>@(dl,ul) \ar@[blue]@/^.5pc/[r]  \ar@[blue]@<0.2mm>@/^.5pc/[r]  & \bullet^{2} \ar@[blue]@(dr,ur) \ar@[blue]@<0.2mm>@(dr,ur) \ar@[blue]@/^.5pc/[l]  \ar@[blue]@<0.2mm>@/^.5pc/[l]  }$$

We have that $B(\a,\b,\g,\d)\cong B(\a',\b',\g',\d')$ if and only if $\exists k_1,k_2\in D^\times$, 
$\a'=k_1\a$, $\b'=\frac{k_1^2}{k_2}\b$, $\g'=\frac{k_2^2}{k_1}\g$, $\d'=k_2\d$. The moduli set is $\bm{\omega}$.
\end{enumerate}
\end{enumerate}
\end{theorem}
From this classification we can recover the classification of two dimensional perfect evolution algebras  over arbitrary fields  in \cite{CGMMS}.
\begin{corollary}\label{oroc}
If $D$ is a field, then the two-dimensional perfect $D$-algebras are: $A_1$, $A_{2,\a}$ (for $\a\ne 0$), $A_{3,\a}$ ($\a\ne 0$), 
$A(\l,1,1)$ ($\l\ne 0$) and $B(1,\l,\m,1)$ ($\l,\m\ne 0,1$, $\l\m\ne 1$).
\end{corollary}

{\bf Data availability.}
The authors confirm that the data supporting the findings of this study are available within the article and its supplementary materials.

\end{document}